\documentclass[12pt]{amsart}
\usepackage{amssymb}
\usepackage{mathrsfs, amssymb, array}
\usepackage{verbatim}
\usepackage[matrix,arrow,curve]{xy}
\usepackage{hhline}
\usepackage{array}

\newcommand{\ZZ}{\mathbb{Z}}
\newcommand{\PP}{\mathbb{P}}
\newcommand{\QQ}{\mathbb{Q}}

\newcommand{\OOO}{\mathscr{O}}
\newcommand{\rk}{\operatorname{rk}}

\newcommand{\Bir}{\operatorname{Bir}}
\newcommand{\Pic}{\operatorname{Pic}}
\newcommand{\Aut}{\operatorname{Aut}}
\newcommand{\Cr}{\operatorname{Cr}}

\newcommand{\xref}[1]{{\rm\ref{#1}}}

 \newcommand{\comp}{\mathbin{\scriptstyle{\circ}}}

\newtheorem{theorem}{Theorem}
\newtheorem{stheorem}[equation]{}

\numberwithin{theorem}{section}
\numberwithin{equation}{theorem}

\newtheorem{mtheorem}[theorem]{}

\theoremstyle{definition}

\newtheorem{say}[theorem]{}

\newtheorem{ppar}[equation]{}

\theoremstyle{definition}

\title{On stable conjugacy of finite subgroups of the plane Cremona group, I}

\author{Fedor Bogomolov}
\address{F. B.: Courant Institute of Mathematical Sciences, 251 Mercer St., New York, NY 10012, U.S.A.}
\email{bogomolo@cims.nyu.edu}
\author{Yuri Prokhorov}
\address{Y. P.:
Steklov Mathematical Institute, 
8 Gubkina str., Moscow 119991, Russia}
\email{prokhoro@gmail.com}
\thanks{
Both authors were partially supported by AG Laboratory SU-HSE, RF~government
grant ag.~11.G34.31.0023.
The second author was partially supported by the grants
RFBR-11-01-00336-a, N.Sh.-5139.2012.1, and
Simons-IUM fellowship.
}

\address{F. B. \& Y. P.: Laboratory of Algebraic Geometry, GU-HSE, 7 Vavilova str.,
Moscow 117312, \mbox{Russia}}

\begin{document}

 \begin{abstract}
We discuss the problem of
 stable conjugacy of finite subgroups of Cremona groups.
We show  that the group $H^1(G,\Pic(X))$
is a stable birational invariant and compute 
this group in some easy cases.
 \end{abstract}
\maketitle
\section{Introduction}
Let $\Bbbk$ be an algebraically closed field of characteristic $0$.
The \textit{Cremona group} $\Cr_n(\Bbbk)$ over $\Bbbk$ is the group of birational
automorphisms of the projective space $\PP^n$, or, equivalently, the group of
$\Bbbk$-automorphisms of the field $\Bbbk(x_1, x_2, \dots, x_n)$ of rational functions
in $n$ independent variables.
Note that embeddings adding new variables one gets a tower $\subset \Cr_{n}(\Bbbk)\subset \Cr_{n+1}(\Bbbk)\subset \cdots$.
Subgroups $G_1,\, G_2\subset \Cr_n(\Bbbk)$ are said to be \textit{stably conjugate} if 
they are conjugate in some $\Cr_{m}(\Bbbk)\supset \Cr_n(\Bbbk)$.
Stable conjugacy of Cremona groups is an analog of 
stable birational equivalence (see \cite{BCTSSD85}, \cite{Voskresenskiui1977}) of 
varieties over non-closed fields.

We are interested in stable conjugacy of subgroups of $\Cr_n(\Bbbk)$.
An example of a subgroup $G\subset \Cr_n(\Bbbk)$ (for $n=3$) that is not stably conjugate to 
any subgroup induced by a linear action was constructed by V. Popov \cite{Popov2011a}.
The construction is based on the example of 
M. Artin and D. Mumford
\cite{Artin-Mumford-1972} and uses non-triviality
of torsions of $H^3(Y,\ZZ)$ for the quotient variety $Y=\PP^3/G$.
In this paper we use another, very simple approach.

Recall that any finite subgroup $G\subset \Cr_n(\Bbbk)$ is induced by
a \emph{biregular} action on a non-singular projective rational variety $X$
(see \ref{G-varieties}).
One can easily show that the group $H^1(G,\Pic(X))$
is stable birational invariant and so it does not depend on the choice of $X$
(see Corollary \ref{Corollary-H1}).
In particular, if $G\subset \Cr_n(\Bbbk)$ is stably conjugate to a linear action, 
then $H^1(G,\Pic(X))=0$.
 This fact is not surprising: in the arithmetic case it was known
for a long time (see e.g. \cite{Manin-1966}, \cite{Voskresenskiui1977}). 
An interesting fact is that in the geometric case the group 
$H^1(G,\Pic(X))$ typically admits a good description.
Note that in  \cite{Lemire-Popov-Reichstein}
the authors used a similar approach to construct non-stably conjugate embeddings 
of certain groups to $\Cr_n(\Bbbk)$, $n\ge 3$
(see \cite[Example 1.36]{Lemire-Popov-Reichstein}).
Our method is  very elementary and can applied to much  wider classes of groups.

We concentrate on the case of the plane Cremona group
$\Cr_2(\Bbbk)$.
In this case the classification of finite subgroups has a long history.
It was started in works of E. Bertini and completed recently
by I. Dolgachev and V. Iskovskikh
\cite{Dolgachev-Iskovskikh} (however even in this case some questions remain open).
We prove the following.

\begin{mtheorem}{\bf Theorem.}\label{theorem-main-p}
Let a finite cyclic group $G$ of prime order $p$ act on a non-singular projective rational surface $X$.
Assume that $G$ fixes \textup(point-wise\textup) a curve of genus $g>0$.
Then
\begin{equation}
\label{main-equality}
H^1(G,\Pic(X))\simeq (\ZZ/p\ZZ)^{2g}.
\end{equation}
\end{mtheorem}

Taking the results of 
\cite{Bayle-Beauville-2000}, \cite{Fernex2004} and \cite{Beauville2004}
into account (see Theorem \ref{Theorem-involutions-Cr2}) we get  the following.

\begin{mtheorem}{\bf Corollary.}\label{Corollary-main-p}
In the notation of \xref{theorem-main-p} the following are 
equivalent:
\begin{enumerate}
\item 
$H^1(G,\Pic(X))=0$,
\item 
$G$ does not fix point-wise a curve $C$ of positive genus,
\item 
$(X,G)$ is conjugate to a linear action on $\PP^2$,
\item 
$(X,G)$ is stably conjugate to a linear action on $\PP^n$ for some $n$.
\end{enumerate}
\end{mtheorem}

In particular, we see that classical de Jonqui\`eres, Bertini and Geiser involutions 
are not stably linearizable.
Another application of the above corollary is that the action 
of the simple Klein group of order 168 on del Pezzo surface of degree $2$ is also 
not stably linearizable.
More application will be given in the forthcoming 
second  part of the paper.
\par\smallskip\noindent
\textbf{Acknowledgments.} This work was completed while the authors were visiting 
University of Edinburgh.
They are grateful to this institution and personally to 
Ivan Cheltsov for invitation and hospitality.
The authors would like to thank A. Beauville, 
J.-L. Colliot-Th{\'e}l{\`e}ne, 
V. L. Popov, and C. Shramov for 
constructive suggestions and criticism.

\section{Preliminaries}

\begin{say}{\bf $G$-varieties.}\label{G-varieties}
Let $G$ be a finite group. A \textit{$G$-variety} is a pair $(X, \rho)$, where $X$ is a
projective variety and $\rho$ is an isomorphism from $G$ to $\Aut(X)$.
A \textit{morphism} (resp. \textit{rational map}) of the pairs 
$(X, \rho) \to (X', \rho')$ is defined to be a
morphism $f : X \to X'$ (resp. rational map  $f : X \dashrightarrow X'$)
such that $\rho' (G' ) = f \comp \rho(G) \comp f^{-1}$.
In particular, two subgroups of $\Aut(X)$
define isomorphic (resp. $G$-birationally equivalent) 
$G$-varieties if and only if they
are conjugate inside of $\Aut(X)$ (resp. $\Bir(X)$). 
If no confusion is likely, we will denote a
$G$-variety by $(X, G)$ or even by $X$.

Now let $X$ be an algebraic variety and let  
$G\subset \Bir(X)$ be a finite subgroup.
By shrinking $X$ we may assume that $G$ acts on $X$ 
by biregular automorphisms. Then replacing $X$ with 
the normalization of $X/G$ in the field $\Bbbk(X)$ we may assume that
$X$ is projective.
Finally we can apply an equivariant resolution of singularities
and replace $X$ with its non-singular  (projective) model.
Thus $G\subset \Bir(X)$ is induced by a biregular action
of $G$ on a non-singular  projective $G$-variety.
In particular, this construction can be applied to 
finite subgroups of $\Cr_n(\Bbbk)$.
\end{say}

\begin{say}{\bf Minimal rational $G$-surfaces.}
Let $X$ be a projective non-singular $G$-surface, where $G$ is a finite group.
It is said to be \textit{$G$-minimal} if any birational
$G$-equivariant 
morphism $f: X\to Y$ is an isomorphism.
It is well-known that for any  $G$-surface there is a minimal 
(projective non-singular) model (see e.g. \cite{Iskovskikh-1979s-e}).
If the surface $X$ is additionally rational, then  
 one of the following holds \cite{Iskovskikh-1979s-e}:
\begin{itemize}
 \item
$X$ is a del Pezzo surface whose invariant Picard number $\Pic(X)^G$ is of rank $1$, or
\item
$X$ admits a structure of $G$-conic bundle, that is, there exists a
surjective $G$-equivariant morphism $f: X\to \PP^1$ such that $f_*\OOO_X=\OOO_{\PP^1}$,
$-K_X$ is $f$-ample and $\rk \Pic(X)^G=2$.
\end{itemize}
\end{say}
\begin{say}{\bf Elements of prime order.}
The classification of elements of prime order in the space Cremona group 
can be summarized as follows.

\begin{stheorem}{\bf Theorem (\cite{Bayle-Beauville-2000}, \cite{Fernex2004}, \cite{Beauville2004}).} 
\label{Theorem-involutions-Cr2}
Let $G=\langle\delta \rangle\subset \Cr_2(\Bbbk)$ is cyclic subgroup of prime order $p$ and let $(X,G)$
be its non-singular projective model. Then the following hold.
\begin{enumerate}
 \item 
The action $(X,G)$ is conjugate to a 
linear action on $\PP^2$ if and only if the fixed point locus  $X^G$ does not contain any curve of positive genus.
 \item 
If $X^G$ contains a curve $C$ of genus $g>0$, then other irreducible components of $X^G$
are either points or rational curves.
In this case the minimal model $(X_{\min},G)$ is unique up to isomorphism and 
there are the following possibilities:
\end{enumerate}
\rm
\renewcommand{\arraystretch}{1.2}
\setlength{\tabcolsep}{8pt}
\begin{center}
\begin{tabular}[]{l|l|l|l|l}
$p$ &$g$&$K_X^2$&$X_{\min}$&$\delta$
\\\hline
$2$& $\ge 1$&$6-2g$&conic bundle&de Jonqui\`eres involution
\\
$2$& $3$&$2$&del Pezzo surface&Geiser involution
\\
$2$& $4$&$1$&del Pezzo surface&Bertini involution
\\
$3$& $1$&$3$&del Pezzo surface&\cite[A1]{Fernex2004}
\\
$3$& $2$&$1$&del Pezzo surface&\cite[A2]{Fernex2004}
\\
$5$& $1$&$1$&del Pezzo surface&\cite[A3]{Fernex2004}
\end{tabular}
\end{center}
\end{stheorem}

\end{say}
\begin{say}{\bf Stable conjugacy.}
Let $(X,G)$ and $(Y,G)$ be $G$-varieties. We say that 
$(X,G)$ and $(Y,G)$ are \textit{stably birational}
if for some $n$ and $m$ there exists an equivariant 
birational map $X\times \PP^n \dashrightarrow Y\times \PP^m$, where 
actions on $\PP^n$ and $\PP^m$ are trivial.
This is equivalent to the conjugacy of the embeddings
$G\subset \Bbbk(X)(x_1,\dots,x_n)$ and  $G\subset \Bbbk(X)(x_1,\dots,x_m)$,
i.e. stable conjugacy of $G\subset \Bbbk(X)$ and  $G\subset \Bbbk(X)$.
\end{say}

The following fact is well-known in the arithmetic case
(see e.g. \cite{Voskresenskiui1977}).
Since we were not able to find a good reference for the present, geometric form,
we provide a complete proof.

\begin{mtheorem}{\bf Proposition (cf. {\cite[2.2]{Manin-1966}}).}
\label{proposition-Pic}
Let $X_1$ and $X_2$ be projective non-singular $G$-varieties.
Assume that $(X_1,G)$ and $(X_2,G)$ are stably birational. Then
there are permutation $G$-modules $\Pi_1$ and $\Pi_2$ such that
the following isomorphism of $G$-modules holds
\[
\Pic(X_1)\oplus \Pi_1\simeq \Pic(X_2)\oplus \Pi_2.
\]
\end{mtheorem}
The proof below is quite standard and depends on the resolution 
of singularities. There is 
more sophisticated  but similar proof due to L. Moret-Bailly
which works in positive characteristic as well
(see \cite[6.2]{Moret-Bailly1986}, \cite[2.A.1]{Colliot-Thelene-Sansuc-1987}).
\begin{proof}
By our assumption, for some $n$ and $m$, there exists a $G$-birational map
$X_1\times \PP^n \dashrightarrow X_2\times \PP^m$, where the
action of $G$ on $\PP^n$ and $\PP^m$ is trivial. Replacing $X_1$ and $X_2$ with
$X_1\times \PP^n$ and $X_2\times \PP^m$ respectively we may assume that
there exists a $G$-birational map $X_1\dashrightarrow X_2$.
Consider a common $G$-equivariant resolution (see e.g. \cite{Abramovich-Wang})
\[
\xymatrix{
&W\ar[rd]^{f_2}\ar[ld]_{f_1}&
\\
X_1\ar@{-->}[rr]&&X_2
}
\]
Then the maps $f_1^*$ and $f_2^*$
induce isomorphisms
$\Pic(W)\simeq \Pic(X_1)\oplus \Pi_1\simeq \Pic(X_2)\oplus \Pi_2$,
where $\Pi_1$ (resp. $\Pi_2$) is a free $\ZZ$-module whose basis is formed by the
prime $f_1$-exceptional (resp. $f_2$-exceptional) divisors.
Since $f_1$ and $f_2$ are $G$-equivariant, the group $G$ permutes these divisors,
so $\Pi_1$ and $\Pi_2$ are permutation modules.
\end{proof}

\begin{stheorem}{\bf Corollary.}\label{Corollary-H1}
In the notation of Proposition \textup{\ref{proposition-Pic}}
we have $H^1(G,\Pic(X_1))\simeq H^1(G,\Pic(X_2))$.
\end{stheorem}
\begin{proof}
By Shapiro's lemma $H^1(G,\Pi_1)= H^1(G,\Pi_2)=0$ (see e.g. \cite[Ch. 1, \S 2.5]{Serre-1962}).
\end{proof}

\begin{stheorem}{\bf Corollary.}
If in the notation of Proposition \textup{\ref{proposition-Pic}}
$(X,G)$ is stably linearizable, then
$H^1(H,\Pic(X))=0$ for any subgroup $H\subset G$.
\end{stheorem}

\section{Proof of Theorem \xref{theorem-main-p}}
Let $\delta\in G$ be a generator and let
$C$ be a (smooth) curve of fixed points with $g:=g(C)>0$.
We  replace $(X,G)$ with its minimal model.
First we consider the case where $X$ is a del Pezzo surface with $\rk\Pic(X)^G=1$.
We start with more general settings.
\begin{say}\label{notation-computation-sequence}
{\bf Del Pezzo case.}
Let $X$ be a del Pezzo surface and let $d:=K_X^2$.
Let $G\subset \Aut(X)$ be any finite subgroup such that
$\Pic(X)^G\simeq \ZZ$. 
Denote 
\[
Q:=\{ x\in \Pic(X) \mid x\cdot K_X=0\}. 
\]

\begin{stheorem}{\bf Lemma.}\label{lemma-computation-sequence}
In the above notation there exists the following
natural exact sequence
\[
0\to \ZZ/d\ZZ \longrightarrow H^1(G,Q) \longrightarrow
H^1(G,\Pic(X)) \longrightarrow 0.
\]
\end{stheorem}
\begin{proof}
By our assumption $g>0$, we have  $d\le 6$ (see e.g. Theorem \ref{Theorem-involutions-Cr2}).
Hence, $\Pic(X)^G = \ZZ\cdot [K_X]$ (see \cite{Iskovskikh-1979s-e}). 
Then the assertion follows from the exact sequence of $G$-modules
\begin{equation}
\label{equation-sequence}
0 \longrightarrow Q \longrightarrow \Pic(X) \overset{\cdot K_X}{\longrightarrow} \ZZ \longrightarrow 0
\end{equation}
because  $H^1(G,\ZZ)=0$
for a finite group $G$.
\end{proof}

\begin{stheorem}{\bf Lemma.}\label{lemma-d=1-computation}
In the notation of \xref{notation-computation-sequence}
assume that  $G$ be a cyclic group generated by $\delta\in G$.
Then the order of $H^1(G,\Pic(X))$ equals to
$|\chi_{\delta,Q}(1)|/d$,
where $\chi_{\delta,Q}(t)$ is the characteristic polynomial of the action of
$\delta$ on $Q$.
\end{stheorem}

The proof uses the following easy observation.
\begin{ppar}{\bf Observation.}\label{observation-notation}
Let $G$ be a cyclic group of order $n$ generated by $\delta\in G$.
Let $\Pi$ be a $\ZZ$-torsion free $\ZZ[G]$-module. Denote
\[
N:=1+\delta+\cdots+\delta^{n-1},\quad
\eta:=1-\delta \in \ZZ[G].
\]
Then 
\[
H^1(G,\Pi)=\ker (N)/ \eta (\Pi).
\]
\end{ppar}
\begin{proof}[Proof of Lemma \textup{\ref{lemma-d=1-computation}}]
Apply the above fact  with $\Pi=Q$. We get $\ker(N)=Q$
(because $N(Q)\subset Q^\delta=0$). Hence
$[Q: \eta(Q)]= |\det (\eta)|$.
\end{proof}

\begin{stheorem}{\bf Corollary.}\label{Corollary-del-Pezzo-computation}
In the notation of \xref{notation-computation-sequence},
let $G$ be a cyclic group of prime order $p$.
Then $d=p^j$, where $j=0$ or $1$, and
\[
H^1(G,\Pic(X))\simeq (\ZZ/p\ZZ)^{(9-d)/(p-1)-j}.
\]
\end{stheorem}
\begin{proof}
Since $Q$ has no $\delta$-invariant vectors, 
the only $\QQ$-irreducible factor of  $\chi_{\delta,Q}(t)$
is the cyclotomic polynomial $t^{p-1}+\cdots+t+1$.
Hence we have
\[
 \chi_{\delta,Q}(t)=(t^{p-1}+\cdots+t+1)^{s},
 \quad 
 s= \rk(Q)/(p-1)=(9-d)/(p-1).
\]
On the other hand, $H^1(G,\Pic(X))$ is a $p$-torsion group
because $G\simeq \ZZ/p\ZZ$
and $\ZZ/d\ZZ$ is a subgroup of 
$\ZZ/p\ZZ$ by \eqref{equation-sequence}.
Hence $d=p^j$ with $j=0$ or $1$ and
$H^1(G,\Pic(X))\simeq (\ZZ/p\ZZ)^{s-j}$.
\end{proof}

By the classification theorem \ref{Theorem-involutions-Cr2}, for $(p,d,g)$
we have one of the following possibilities:
$(2,1,4)$, $(2,2,3)$, 
$(3,3,1)$, $(3,1,2)$, $(5,1,1)$.
Applying Corollary \ref{Corollary-del-Pezzo-computation}
we get our equality \eqref{main-equality}.
\end{say}

\begin{say}{\bf Conic bundle case.}
Now assume that $X$ has a structure of $G$-equivariant conic bundle $f: X\to\PP^1$.
Then again by Theorem \ref{Theorem-involutions-Cr2}\
$\delta$ is a \emph{de Jonqui\`eres involution} of genus $g>0$.
Recall that this is an element
$\delta\in \Cr_2(\Bbbk)$ of order $2$ induced by an action on
a ($G$-equivariant) relatively minimal conic bundle $f: X\to \PP^1$
with $2g+2$ degenerate fibers so that the locus of fixed points is a hyperelliptic curve $C$ of genus $g$
(elliptic curve if $g=1$) and the restriction $f|_C: C\to \PP^1$
is a double cover (see \cite{Bayle-Beauville-2000},
\cite{Dolgachev-Iskovskikh} for details).
\end{say}

 \begin{say} \label{de-Jonquieres-involutions}
Let $F$ be a typical fiber and let
$F_i=F_i'+F_i''$, $i=1,\dots, 2g+2$ be all the degenerate fibers.
Let $Q\subset \Pic(X)$ be the $\ZZ$-submodule of rank $2g+3$ generated by the 
components of degenerate fibers.
It has a $\ZZ$-basis consisting of
$F$ and $F_i'$, $i=1,\dots, 2g+2$.
The action of $\delta$ is given by
\[
\delta : F_i'\longmapsto F-F_i',\quad F \longmapsto F.
\]
Apply  \ref{observation-notation} with $\Pi=Q$. We have
\[
\ker(N)= \left\{\alpha F +\sum \alpha_i F_i'\in Q
\ \left | \ 2\alpha +\sum \alpha_i=0 \right.\right\}
\]
and $\eta(Q)$ is generated by the classes of $2F_i'-F$.
Therefore,
\[
H^1(G,Q)=(\ZZ/2\ZZ)^{2g+1}.
\]
Since $H^1(G,\ZZ)=0$, from the exact sequence
\begin{equation}\label{exact-sequence-conic-bundle}
0\longrightarrow Q \longrightarrow\Pic(X) \overset{\cdot F}\longrightarrow \ZZ \longrightarrow 0
\end{equation}
we get
\[
0\to Q^G \longrightarrow\Pic(X)^G \overset{\cdot F}\longrightarrow \ZZ
\longrightarrow H^1(G,Q) \longrightarrow H^1(G,\Pic(X)) \to 0.
\]
Note that the image of $\Pic(X)^G$ in $\ZZ$ is generated by $K_X\cdot F=2$.
Therefore, $H^1(G,\Pic(X))=(\ZZ/2\ZZ)^{2g}$.
 \end{say}

This proves Theorem \xref{theorem-main-p}.

\begin{proof}[Proof of Corollary \xref{Corollary-main-p}.]
The implication (i) $\Rightarrow$ (ii) follows by Theorem \ref{theorem-main-p},
implications (iii) $\Rightarrow$ (iv) $\Rightarrow$ (i) are  obvious, and
 (ii) $\Rightarrow$ (iii) follows by Theorem \ref{Theorem-involutions-Cr2}.
\end{proof}

\def\cprime{$'$} \def\polhk#1{\setbox0=\hbox{#1}{\ooalign{\hidewidth
  \lower1.5ex\hbox{`}\hidewidth\crcr\unhbox0}}}


\begin{thebibliography}{BCTSSD85}

\bibitem[AW97]{Abramovich-Wang}
D. Abramovich and J. Wang.
\newblock {Equivariant resolution of singularities in characteristic 0.}
\newblock {\em Math. Res. Lett.}, 4(2-3):427--433, 1997.

\bibitem[AM72]{Artin-Mumford-1972}
M.~Artin and D.~Mumford.
\newblock Some elementary examples of unirational varieties which are not
  rational.
\newblock {\em Proc. London Math. Soc. (3)}, 25:75--95, 1972.

\bibitem[BB00]{Bayle-Beauville-2000}
L. Bayle and A. Beauville.
\newblock Birational involutions of {${\bf P}\sp 2$}.
\newblock {\em Asian J. Math.}, 4(1):11--17, 2000.
\newblock Kodaira's issue.

\bibitem[BB04]{Beauville2004}
A. Beauville and J. Blanc.
\newblock On {C}remona transformations of prime order.
\newblock {\em C. R. Math. Acad. Sci. Paris}, 339(4):257--259, 2004.

\bibitem[BCSS85]{BCTSSD85}
A. Beauville, J.-L. Colliot-Th{\'e}l{\`e}ne, J.-J. Sansuc, and
  P. Swinnerton-Dyer.
\newblock Vari\'et\'es stablement rationnelles non rationnelles.
\newblock {\em Ann. of Math. (2)}, 121(2):283--318, 1985.


\bibitem[CTS87]{Colliot-Thelene-Sansuc-1987}
J.-L. Colliot-Th{\'e}l{\`e}ne and J.-J. Sansuc.
\newblock La descente sur les vari\'et\'es rationnelles. {II}.
\newblock {\em Duke Math. J.}, 54(2):375--492, 1987.

\bibitem[dF04]{Fernex2004}
T. de~Fernex.
\newblock On planar {C}remona maps of prime order.
\newblock {\em Nagoya Math. J.}, 174:1--28, 2004.

\bibitem[DI09]{Dolgachev-Iskovskikh}
I.~V. Dolgachev and V.~A. Iskovskikh.
\newblock Finite subgroups of the plane {C}remona group.
\newblock In {\em Algebra, arithmetic, and geometry: in honor of {Y}u. {I}.
  {M}anin. {V}ol. {I}},  {\em Progr. Math.}, v. 269, 443--548.
  Birkh\"auser Boston Inc., Boston, MA, 2009.

\bibitem[Isk80]{Iskovskikh-1979s-e}
V.~A. Iskovskikh.
\newblock Minimal models of rational surfaces over arbitrary fields.
\newblock {\em Math. USSR-Izv.}, 14(1):17--39, 1980.

\bibitem[LPR06]{Lemire-Popov-Reichstein}
N. Lemire, V.~L. Popov, and Z. Reichstein.
\newblock Cayley groups.
\newblock {\em J. Amer. Math. Soc.}, 19(4):921--967, 2006.

\bibitem[Man66]{Manin-1966}
Yu.~I. Manin.
\newblock Rational surfaces over perfect fields.
\newblock {\em Inst. Hautes \'Etudes Sci. Publ. Math.}, (30):55--113, 1966.

\bibitem[MB86]{Moret-Bailly1986}
L. Moret-Bailly.
\newblock Vari\'et\'es stablement rationnelles non rationnelles (d'apr\`es
  {B}eauville, {C}olliot-{T}h\'el\`ene, {S}ansuc et {S}winnerton-{D}yer).
\newblock {\em Ast\'erisque}, (133-134):223--236, 1986.
\newblock Seminar Bourbaki, Vol. 1984/85.


\bibitem[Pop11]{Popov2011a}
V.~L. Popov.
\newblock Some subgroups of the {C}remona groups.
\newblock {\em arXiv e-print}, arXiv:1110.2410, 2011.

\bibitem[Ser63]{Serre-1962}
J.-P. Serre.
\newblock {\em Cohomologie galoisienne},  {\em Cours au Coll\`ege
  de France}, v. 1962.
\newblock Springer-Verlag, Berlin, 1962/1963.

\bibitem[Vos77]{Voskresenskiui1977}
V.~E. Voskresenski{\u\i}.
\newblock {\em Algebraicheskie tory}.
\newblock Izdat. ``Nauka'', Moscow, 1977.

\end{thebibliography}
\end{document}